%% The Degenerate (Hyperbolic) Switch 
%update 26 September 2006 previous to old
\magnification=\magstep1

\def\nl{\hfill\break}

\def\hexnumber#1{\ifcase#1 0\or 1\or 2\or 3\or 4\or 5\or 6\or 7\or 8\or
 9\or A\or B\or C\or D\or E\or F\fi}
%
%     Next we define the AMS symbol-a fonts at 12,10,9,7,6,5 points
%
\font\twelvemsa=msam10 scaled 1200   %%%  ** change to 1315 if necessary
\font\tenmsa=msam10                  %%%  see note above
\font\ninemsa=msam9            \font\sevenmsa=msam7
\font\sixmsa=msam6             \font\fivemsa=msam5
%%%%%%  (add further sizes here if you need them)
%
%    and the standard size family for these fonts
%
\newfam\msafam                 \textfont\msafam=\tenmsa
\scriptfont\msafam=\sevenmsa   \scriptscriptfont\msafam=\fivemsa
\edef\hexa{\hexnumber\msafam}        %  The msa family is  \fam\hexa
\def\msa{\fam\msafam\tenmsa}         %  \msa  switches to this family
%
%    Repeat these steps for the AMS symbol-b fonts
%
\font\twelvemsb=msbm10 scaled 1200   %%%  ** change to 1315 if necessary
\font\tenmsb=msbm10                  %%%  see note above
\font\ninemsb=msbm9            \font\sevenmsb=msbm7
\font\sixmsb=msbm6             \font\fivemsb=msbm5
%%%%%%  (add further sizes here if you need them)
%
\newfam\msbfam                 \textfont\msbfam=\tenmsb       
\scriptfont\msbfam=\sevenmsb   \scriptscriptfont\msbfam=\fivemsb
\edef\hexb{\hexnumber\msbfam}        %  The msb family is \fam\hexb 
\def\msb{\fam\msbfam\tenmsb}         %  \msb switches to this family
%
%        Repeat for the Euler-Fraktal fonts 
%
\font\twelveeufm=eufm10 scaled 1200  %%%  ** change to 1315 if necessary
\font\teneufm=eufm10                 %%%  see note above
\font\nineeufm=eufm9           \font\seveneufm=eufm7
\font\sixeufm=eufm6            \font\fiveeufm=eufm5
%%%%%%  (add further sizes here if you need them)
%
\newfam\eufmfam                \textfont\eufmfam=\teneufm
\scriptfont\eufmfam=\seveneufm \scriptscriptfont\eufmfam=\fiveeufm
\edef\hexf{\hexnumber\eufmfam}      % The Euler-Fraktal family is
\def\frak{\fam\eufmfam\teneufm}     % \fam\hexf and \frak switches to this
%
%%%  Add further fonts families here (using the same format) if you need
%    them.  The def of hexnumber is optional (it is only used for
%    \mathchardef 's).
% 
%       Now we need to define the standard fonts (which are
%       already defined at 10,7 and 5 point) at 12,9 and 6 point:
%
%      Roman fonts:
\font\twelverm=cmr10 scaled 1200    %%%  ** change to 1315 if necessary
\font\ninerm=cmr9                   %%%  see note above
\font\sixrm=cmr6   
%%%%%%  (add further sizes here if you need them)
%
%      Math italic fonts
\font\twelvei=cmmi10 scaled 1200    %%%  ** change to 1315 if necessary
\font\ninei=cmmi9                   %%%  see note above
\font\sixi=cmmi6  
%%%%%%  (add further sizes here if you need them)
%
%     Symbol fonts
\font\twelvesy=cmsy10 scaled 1200   %%%  ** change to 1315 if necessary
\font\ninesy=cmsy9                  %%%  see note above
\font\sixsy=cmsy6  
%%%%%%  (add further sizes here if you need them)
%
%     Bold face
\font\twelvebf=cmbx10 scaled 1200   %%%  ** change to 1315 if necessary
\font\ninebf=cmbx9                  %%%  see note above
\font\sixbf=cmbx6  
%%%%%%  (add further sizes here if you need them)
%
%     Finally three fonts (text italic, slanted and typewriter type)
%     which are only defined at text size
%
\font\twelveit=cmti10 scaled 1200   %%%  ** change to 1315 if necessary
\font\nineit=cmti9                  %%%  see note above
%%%%%%  (add further sizes here if you need them)
%
\font\twelvesl=cmsl10 scaled 1200   %%%  ** change to 1315 if necessary
\font\ninesl=cmsl9                  %%%  see note above
%%%%%%  (add further sizes here if you need them)
%
\font\twelvett=cmtt10 scaled 1200   %%%  ** change to 1315 if necessary
\font\ninett=cmtt9                  %%%  see note above
%%%%%%  (add further sizes here if you need them)
%
%
%     Now come the two main macros.  What  \small  does is to
%     change all the families of fonts from normal size which is
%     10,7,5  (ie 10pt text, 7pt subscript, 5pt subsubscript)
%     to 9,6,5.  \large  similarly changes to  12,9,7.  To make
%     other size changing macros, choose your three sizes, add
%     font size definitions if necessary and make the obvious changes
%     to one of these macros.  Change  \normalbaselineskip  and
%     \strutbox  dimensions to appropriate sizes as well.  To
%     add further fonts, insert them in each macro, using the
%     AMS fonts as a model.
%      
%
\def\small{%
%
%   redefine the sizes of the roman fonts :
%
\textfont0=\ninerm \scriptfont0=\sixrm \scriptscriptfont0=\fiverm
\def\rm{\fam0\ninerm}        % ( \rm  sets \ninerm  in text mode
%                            %  and \fam0 in math mode)
%
%   and the math italic fonts :
%
\textfont1=\ninei \scriptfont1=\sixi \scriptscriptfont1=\fivei
%
%   and the symbol fonts :
%
\textfont2=\ninesy \scriptfont2=\sixsy \scriptscriptfont2=\fivesy
%
%   There is only one math extension font :
%
\textfont3=\tenex \scriptfont3=\tenex \scriptscriptfont3=\tenex
%
%   Next the bold font (named rather than numbered) :
%
\textfont\bffam=\ninebf \scriptfont\bffam=\sixbf
\scriptscriptfont\bffam=\fivebf \def\bf{\fam\bffam\ninebf}%
%
%   and the three text-only fonts : 
%
\textfont\itfam=\nineit \def\it{\fam\itfam\nineit}%
\textfont\slfam=\ninesl \def\sl{\fam\slfam\ninesl}%
\textfont\ttfam=\ninett \def\tt{\fam\ttfam\ninett}%
%
%   Now the three new families of AMS fonts :
%
%   AMS symbol-a
%
\textfont\msafam=\ninemsa \scriptfont\msafam=\sixmsa
\scriptscriptfont\msafam=\fivemsa \def\msa{\fam\msafam\ninemsa}%         
%
%   AMS symbol-b
%
\textfont\msbfam=\ninemsb \scriptfont\msbfam=\sixmsb
\scriptscriptfont\msbfam=\fivemsb \def\msb{\fam\msbfam\ninemsb}%         
%
%   Euler-Fraktal font
%
\textfont\eufmfam=\nineeufm  \scriptfont\eufmfam=\sixeufm
\scriptscriptfont\eufmfam=\fiveeufm \def\frak{\fam\eufmfam\nineeufm}%
%
%%%  Add further fonts families here if you need them.
%
%   Finally reset \normalbaselineskip and \strubox :
%
\normalbaselineskip=11pt
\setbox\strutbox=\hbox{\vrule height8pt depth3pt width0pt}%
%
%   and finish by setting \normalbaselines and \rm (roman) as defaults :
%
\normalbaselines\rm}    %%%   End of  \small  macro      
%
%
%    The \large  macro is similar (comments abbreviated):
%
%
\def\large{%
\textfont0=\twelverm \scriptfont0=\ninerm \scriptscriptfont0=\sevenrm
\def\rm{\fam0\twelverm}%
\textfont1=\twelvei \scriptfont1=\ninei \scriptscriptfont1=\seveni
\textfont2=\twelvesy \scriptfont2=\ninesy \scriptscriptfont2=\sevensy
\textfont3=\tenex \scriptfont3=\tenex \scriptscriptfont3=\tenex
\textfont\bffam=\twelvebf \scriptfont\bffam=\ninebf
\scriptscriptfont\bffam=\sevenbf \def\bf{\fam\bffam\twelvebf}%
\textfont\itfam=\twelveit \def\it{\fam\itfam\twelveit}%
\textfont\slfam=\twelvesl \def\sl{\fam\slfam\twelvesl}%
\textfont\ttfam=\twelvett \def\tt{\fam\ttfam\twelvett}%
%   AMS symbol-a  :
\textfont\msafam=\twelvemsa \scriptfont\msafam=\ninemsa
\scriptscriptfont\msafam=\sevenmsa \def\msa{\fam\msafam\twelvemsa}         
%   AMS symbol-b  :
\textfont\msbfam=\twelvemsb \scriptfont\msbfam=\ninemsb
\scriptscriptfont\msbfam=\sevenmsb \def\msb{\fam\msbfam\twelvemsb}         
%   Euler-Fraktal font :
\textfont\eufmfam=\twelveeufm  \scriptfont\eufmfam=\nineeufm
\scriptscriptfont\eufmfam=\seveneufm \def\frak{\fam\eufmfam\teneufm}
%%%% Add further fonts families here if you need them.
%%   Finally reset \normalbaselineskip and \strubox and initialise :
\normalbaselineskip=15pt
\setbox\strutbox=\hbox{\vrule height11pt depth4pt width0pt}%
\normalbaselines\rm}%
%     
%   The next two lines define commonly used switches for
%   blackboard bold (\Bbb) and gothic type (\goth)  
\def\Bbb{\msb}

%
%    To use the new AMS fonts you can either use the control
%    sequences \msa \msb (alias \Bbb) \frak (alias \goth)  eg :

   % see the msam font table
%
%   or, more generally, make \mathchardef's (cf Knuth p155) eg :
\mathchardef\plussquare="0\hexa01
\mathchardef\nge="3\hexb0B
\mathchardef\maltesecross="0\hexa7A
\mathchardef\del="0\hexf01
%
%   Finally, for continuity with earlier definitions, we define 
%   \ninepoint  as an alias for  \small :

%

%\input epsf
%\nopagenumbers
\overfullrule=0pt

\font\Bbb=msbm10

%\font\secfont=cmbx10 scaled \magstep 2
\font\secfont=cmbx10

\font\nam=cmr8
\font\aff=cmti8
\font\refe=cmr9

\mathchardef\square="0\hexa03
\def\qed{\hfill$\square$\par\rm}

\def\boxing#1{\ \lower 3.5pt\vbox{\vskip 3.5pt\hrule \hbox{\strut\vrule
\ #1 \vrule} \hrule} }

\def\down#1{\ \lower 3.5pt\vbox{\vskip 3.5pt \hbox{\strut \ #1 \vrule} \hrule} }
\def\negdown#1{\ \lower 3.5pt\vbox{\vskip 3.5pt \hbox{\strut  \vrule \ #1 }\hrule} }
\def\adj{\hbox{\rm adj}}
\def\tr{\hbox{\rm tr}}

\hsize=6.3 truein
\vsize=9 truein

\baselineskip=13 pt
\parskip=\baselineskip
 1

\parindent=0pt

\def\H{\hbox{\Bbb H}}
\def\E{\hbox{\Bbb E}}

\def\a{\hbox{\bf a}}
\def\b{\hbox{\bf b}}
\def\c{\hbox{\bf c}}
\def\x{\hbox{\bf x}}

\def\i{\hbox{\bf i}}
\def\j{\hbox{\bf j}}
\def\k{\hbox{\bf k}}

%\font\chapfont=cmbx10 scaled \magstep 2

%\font\secfont=cmbx10 scaled \magstep 1

%\hsize=6.2 in

%\vsize=8 in

%\voffset=.75 in

%\parskip=\baselineskip

\newif \iftitlepage \titlepagetrue
%\headline={\iftitlepage   
%\global\titlepagefalse\else \centerline{ Research  Report}\fi}

\def\diagram{\global\advance\diagramnumber by 1
$$\epsfbox{stevefig.\number\diagramnumber}$$}
\def\ddiagram{\global\advance\diagramnumber by 1
\epsfbox{stevefig.\number\diagramnumber}}

\newcount\diagramnumber
\diagramnumber=0

\newcount\secnum \secnum=0
\newcount\subsecnum
\newcount\defnum
\def\section#1{
                \vskip 10 pt
                \advance\secnum by 1 \subsecnum=0
                \leftline{\secfont \the\secnum \rm\quad#1}
                }

\def\subsection#1{
                \vskip 10 pt
                \advance\subsecnum by 1 
                \defnum=1
                \leftline{\secfont \the\secnum.\the\subsecnum\ \rm\quad #1}
                }

\def\definition{
                \advance\defnum by 1 
                \bf Definition 
\the\secnum .\the\defnum \rm \ 
                }

\def\lemma{
                \advance\defnum by 1 
                \par\bf Lemma  \the\secnum
.\the\defnum \rm \ \par
                }

\def\theorem{
                \advance\defnum by 1 
                \par\bf Theorem  \the\secnum
.\the\defnum \rm \ 
               }

\def\cite#1{
				\secfont [#1]
				\rm
}

\vglue 20 pt

\centerline{\secfont QUATERNION ALGEBRAS and}
\centerline{\secfont INVARIANTS of VIRTUAL KNOTS and LINKS}
\centerline{\secfont II: The Hyperbolic Case}

\medskip

\centerline{\nam STEPHEN BUDDEN, ROGER FENN${}^1$}
\centerline{\aff ${}^1$School of Mathematical Sciences, University of Sussex}
\centerline{\aff Falmer, Brighton, BN1 9RH, England}
\centerline{\aff e-mail addresses: rogerf@sussex.ac.uk, stevie\_hair@yahoo.com}
\baselineskip=10 pt
\parskip=0 pt
\bigskip
\centerline{\nam ABSTRACT}
\leftskip=0.25 in
\rightskip=0.25in

{\nam Let $A,\ B$ be invertible, non-commuting elements of a ring
$R$. Suppose that $A-1$ is also invertible and that the equation
$$[B,(A-1)(A,B)]=0$$ called the fundamental equation is
satisfied. Then an invariant $R$-module is defined for any diagram of
a (virtual) knot or link. Solutions in the classic quaternion case
have been found by Bartholomew, Budden and Fenn. Solutions in the
generalised quaternion case have been found by Fenn in an earlier paper.
These latter
solutions are only partial in the case of $2\times2$ matrices and the
aim of this paper is to provide solutions to the missing cases.}

\leftskip=0 in
\rightskip=0 in
\baselineskip=13 pt
\parskip=\baselineskip

\parskip=\baselineskip
\section{Introduction}

 Let $A,\ B$ be invertible, non-commuting elements of a ring
$R$. Suppose that $A-1$ is also invertible. Our aim is to find solutions to 
the equation
$$A^{-1}B^{-1}AB-B^{-1}AB=BA^{-1}B^{-1}A-A,$$
called the {\sl fundamental equation}.  In this case an invariant
$R$-module is defined for any diagram of a (virtual) knot or
link. Solutions in the classic quaternion case have been found, see
\cite{BF}, \cite{BuF}. Solutions in the generalised quaternion case have been
found, see \cite{F}. These are only partial in the case of $2\times2$
matrices; that is sufficient but not necessary conditions were
given. The matrices satisfying these conditions were called {\sl
matching}.  The aim of this paper is to provide solutions of the
missing matrices: the mismatching or hyperbolic matrices.  This means
that this paper together with earlier papers provides all $2\times 2$
matrix solutions to the fundamental equation. These solutions may be
summed up with the help of the following theorem.
\theorem{Suppose $A, B$ are two non-commuting $2\times2$ matrix solutions
of the fundamental equation. Then either
$$\tr(A)=\det(A)\hbox{ and }\tr(AB^{-1})=0$$
or the pair $A, B$ are similar to a pair of the form
 $$A = \pmatrix{a_0 + a_3  &  2a_1\cr  0 & a_0 - a_3 \cr}\qquad
 B =  \pmatrix{ {2b_3\over a_0-a_3}   & 2b_1 \cr 0 &  2b_3
\bigl({1\over a_0-a_3}-2\bigr)\cr}
$$ 
There are a number of sporadic $n\times n$ matrix
solutions to the fundamental equation which will appear in a further paper
with V. Turaev.  However these are almost
certainly not complete and finding the general $n \times n$ solution
is probably very hard.

With these solutions whole new families of invariant modules and
polynomials of virtual knots and links are defined. An appendix
where these are calculated from the examples in the table of N. Kamada
will be put on the web.

Many of the conventions and notation can be found in \cite{F}. We will
reproduce the details necessary for the understanding of this paper and
leave fine details for the interested reader in \cite{F}.

\section{Generalised Quaternions}
Let $F$ be a field of characteristic not equal to 2. Pick two non-zero
elements $\lambda,\ \mu$ in $F$. Let
$\left({{\lambda,\ \mu}\over F}\right)$
denote the algebra of dimension 4 over $F$ with basis $\{1,\i,\j,\k\}$
and relations $\i^2=\lambda,\ \j^2=\mu,\ \i\j=-\j\i=\k$.  The multiplication
table is given by
$$\bordermatrix{&\i&\j&\k\cr \i&\lambda&\k&\lambda \j\cr \j&-\k&\mu&-\mu \i\cr
\k&-\lambda \j&\mu \i&-\lambda\mu\cr}.$$
We will only consider the case $\lambda=-1$ and $\mu=1$. In this case the
algebra is the ring of $2\times 2$
matrices, $M_2(F)=\left({{-1,\ 1}\over F}\right)$ or $\left({{1,\
1}\over F}\right)$. This is the only quaternion algebra with zero
divisors or non-zero isotropic elements.

The generators of $\left({{-1,\ 1}\over F}\right)$ are, together with
the identity, the Pauli matrices
$$\i=\pmatrix{0&1\cr -1&0\cr},\ \j=\pmatrix{0&1\cr 1&0\cr},\ 
\k=\pmatrix{1&0\cr 0&-1\cr}.$$
By an abuse of notation we will often confuse the scalar matrix
$\pmatrix{\nu&0\cr 0&\nu\cr}$ with the corresponding field element
$\nu$.

We denote general $2\times2$ matrices by capital roman letters such as $A, B,
\ldots$. If the trace of a matrix is zero then we denote it by bold face
lower case, $\a,\b,\ldots$ etc.  Field
elements, (scalars) will be denoted by lower case roman letters such
as $a, b, \ldots$ and lower case greek letters such as $\alpha, \beta,
\ldots$. Therefore any $2\times2$ matrix can be written uniquely as
$$A=a_0+\a.$$
We call $\a$ the {\bf traceless part} of $A$. Any traceless matrix can
be written as a unique three dimensional linear combination of the Pauli
matrices.

The {\bf conjugate} of $A$ (sometimes called the adjugate) is
$\overline{A}=a_0-\a$. In symbols
$$\overline{\pmatrix{a&b\cr c&d\cr}}=\pmatrix{d&-b\cr -c&a\cr}$$
The {\bf determinant} of $A$ is
$\det(A)=A\overline{A}$ and the {\bf trace} of $A$ is
$\tr(A)=A+\overline{A}$.

Conjugation is an anti-isomorphism of order 2. That is it satisfies
$$\overline{A+B}=\overline{A}+\overline{B},\quad
\overline{AB}=\overline{B}\;\overline{A},
\quad\overline{aA}=a\overline{A},\quad \overline{\overline{A}}=A.$$
Also $\overline{A}=A$ if and only if $A$ is a scalar and
$\overline{A}=-A$ if and only if $A$ has trace zero.

The determinant is a scalar satisfying $\det(AB)=\det(A)\det(B)$. We
will denote the set of values of the determinant function by ${\cal
N}$. It is a multiplicatively closed subset of $F$ and ${\cal
N}^*={\cal N}-\{0\}$ is a multiplicative subgroup of $F^*$. An element
$A$ has an inverse if and only if $\det(A)\ne 0$ in which case
$A^{-1}=\det(A)^{-1}\overline{A}$.

The trace of a $2\times2$ matrix is twice its scalar part.
A general matrix can be written uniquely as
$$\pmatrix{\alpha&\beta\cr \gamma&\delta\cr}={1\over2}\left[(\alpha+\delta)+
(\beta-\gamma)\i
+(\beta+\gamma)\j+(\alpha-\delta)\k\right]$$
Conversely
$$A=a_0+a_1\i+a_2\j+a_3\k=\pmatrix{a_0+a_3&a_2+a_1\cr a_2-a_1&a_0-a_3\cr}$$

Conjugation is
$$\overline{A}=\adj A=\pmatrix{\delta&-\beta\cr -\gamma&\alpha\cr}=
\pmatrix{a_0-a_3&-a_2-a_1\cr a_1-a_2&a_0+a_3\cr}$$
and the determinant is
$$\det A=\alpha\delta-\beta\gamma=a_0^2+a_1^2-a_2^2-a_3^2$$

The scalar part of $A$ is $a_0=\tr A/2=(\alpha+\delta)/2$ and the traceless
part is
$$\pmatrix{a_3&a_2+a_1\cr a_2-a_1&-a_3\cr}=\pmatrix{(\alpha-\delta)/2&\beta\cr 
\gamma&(\delta-\alpha)/2\cr}$$

\subsection{Various Multiplications on $2\times2$ Matrices}

Let $A, B$ be two $2\times2$ matrices .There is a bilinear form given by 
$$A\cdot B={1\over 2}(A\overline{B}+B\overline{A})={1\over 2}
(\overline{A}B+\overline{B}A)={1\over 2}\tr(A\overline{B}).$$
The corresponding quadratic form is
$\det(A)$.
Let $\a, \b$ be traceless $2\times2$ matrices. Then
$$\a\b=-\a\cdot\b+\a\times\b$$
where $\a\cdot\b$ is
the restriction of the bilinear form and
$\a\times\b$ is the {\bf cross product}.
The cross product has the usual rules of bilinearity and skew symmetry.
The triple cross product expansion
$$\a\times(\b\times\c)=(\c\cdot\a)\b-(\b\cdot\a)\c$$
is easily verified.
The {\bf scalar triple product} is 
$$[\a,\b,\c]=\a\cdot(\b\times\c)=-
\Biggl|
\matrix{a_1&a_2&a_3\cr b_1&b_2&b_3\cr c_1&c_2&c_3\cr}\Biggr|$$
from which all the usual rules can be deduced.
Here $\a=a_1\i+a_2\j+ a_3\k$ etc.

If
$$A=\pmatrix{\alpha_1&\alpha_2\cr 
\alpha_3&\alpha_4\cr}\hbox{ and }B=\pmatrix{\beta_1&\beta_2\cr
\beta_3&\beta_4\cr}\hbox{ then }
A\cdot B={1\over2}(\alpha_1\beta_4-\alpha_2\beta_3-\alpha_3\beta_2+
\alpha_4\beta_1)$$
If $\a=\pmatrix{\alpha_1&\alpha_2\cr \alpha_3&-\alpha_1\cr}$ and
$\b=\pmatrix{\beta_1&\beta_2\cr
\beta_3&-\beta_1\cr}$ are traceless then
$$\a\cdot\b=-\alpha_1\beta_1-(\alpha_2\beta_3+\alpha_3\beta_2)/2$$
and
$$\a\times\b=\pmatrix{(\alpha_2\beta_3-\alpha_3\beta_2)/2&\alpha_1\beta_2-
\alpha_2\beta_1
\cr\alpha_3\beta_1-\alpha_1\beta_3 &(\alpha_3\beta_2-\alpha_2\beta_3)/2\cr}$$

\subsection{Dependancy Criteria}

In this subsection we will consider conditions for sets of $2\times2$
matrices to be linearly dependant or otherwise.  A non-zero $2\times2$
matrix, $A$, is called {\bf isotropic} or {\bf degenerate} if
$\det(A)=0$ and {\bf anisotropic} otherwise. So only non-zero
anisotropic matrice have inverses.

\lemma{A pair of traceless matrices $\a, \b$ is linearly dependant if and
only if $\a\times \b=0$.}

{\bf Proof} The proof is clear one way using the antisymmetry of the
cross product. Conversely suppose $\a\times \b=0$.  Then $(\a\times
\b)\times\c=(\a\cdot\c)\b-(\b\cdot\c)\a=0$. This can be made into a
linear dependancy by a suitable choice of $\c$, for example if
$\a\cdot\c\ne0$.
\qed
As a corollary we have the following
\lemma{Two $2\times2$ matrices  commute if and only their traceless parts
are linearly dependant.} \qed
Now we look for conditions for the triple of traceless matrices,
$\a,\b,\a\times \b$, to be linearly dependant.  The
required condition is given by the following lemma.

\lemma{The traceless matrices $\a,\b,\a\times \b$, are linearly dependant if
and only if 
$$\det(\a)\det(\b)=(\a\cdot\b)^2.$$ This is equivalent to the condition that
$\a\times\b$ is isotropic or zero.}

{\bf Proof}
Three 3-dimensional vectors are linearly dependant if and only if the
determinant they form by rows is zero. In the case of traceless matrices
this means the scalar triple product is zero
$$[\a,\b,\c]=\a\cdot(\b\times\c)=0.$$
Replacing $\c$ with $\a\times\b$ and expanding out using the triple
cross product formula gives the first equation. Using the expansion formul\ae\ 
$$\det(\a\times\b)=\det(\a)\det(\b)-(\a\cdot\b)^2$$ 
gives the second condition.
\qed
We have the following corollary.

\lemma{If $\a, \b$ are traceless matrices and
$\a\times \b$ has non-zero determinant, then the triple
$\a,\b,\a\times \b$, is linearly independant.}
\qed
\section{Solving the Fundamental Equation}
Given a set $X$ let $S$ be an endomorphism of $X^2$. 
Such an $S$ is
called a {\bf switch} if
{\parindent=20pt
\item{1} $S$ is invertible and 
\item{2} the set theoretic Yang-Baxter equation
$$ (S\times id)(id\times
S)(S\times id)=(id\times S)(S\times id)(id\times S)$$
is satisfied. }
Switches are used in \cite{FJK} to define biracks and biquandles by the formula
$$S(a,b)=(b_a,a^b).$$
Switches can be used to find representations of the virtual braid groups and
invariants of virtual knots and links, see \cite{FJK},\cite{BF},\cite{BuF}
and \cite{F}.

We are looking for linear solutions. That is $S=\pmatrix{A & B \cr C &
D\cr}$, where the matrix entries $A, B, C, D$ are elements of $R$, an
associative but not necessarily commutative ring and $X$ is a left
$R$-module.

The solutions can be divided into two types when the entries are not
zero divisors.

{\bf The commutative case}
$$\leqalignno{
\hbox{ The identity }&&0:\cr
&&\cr
S=\pmatrix{0 & B \cr C & 1-BC\cr}& \hbox{ or }
S=\pmatrix{1-BC & B \cr C & 0\cr}&1:\cr}$$
where $B$ and $C$ are arbitary commuting invertible elements.

{\bf The non-commutative case}
$$\leqalignno{
S=\pmatrix{A & B \cr C & D\cr}&&2:\cr}$$ where
$A, A-1, B$ are invertible, $A, B$ do not commute and satisfy the
fundamental equation
$$A^{-1}B^{-1}AB-B^{-1}AB=BA^{-1}B^{-1}A-A$$
moreover
$$C=A^{-1}B^{-1}A(1-A),\ D=1-A^{-1}B^{-1}AB.$$
There are also similar solutions where $A, D$ and $B, C$ are interchanged.
We are only interested in this last case
and are therefore looking for $2\times2$ matrices
$A=a_0 + \a $ and $B=b_0 + \b$ which satisfy the fundamental equation.
Since  $A$ and $B$ do not commute $\a $ and $\b$
are linearly independent.  As in \cite{BuF} and \cite{F} the linear relation
$$\eqalignno{
&(\tr(A)-\det(A))\det(\b)\a+(\det(A)-\tr(A))(\a\cdot\b)\b\cr
&\qquad +(b_0(\det(A)-\tr(A))+2A\cdot B)\a\times\b=0 &(1)\cr}$$
holds.
          
The paper \cite{F} has solved the {\bf matching solutions}. That is
solutions where $\det(A)=\tr(A)$ and $A \cdot B=0$.  So we are interested
in the {\bf mismatching solutions}. In this case $A, B$ satisfy the
fundamental equation and $\a $, $\b $ and $\a \times \b $ are linearly
dependent.

\section{Finding Linearly Dependent Triples  $\a ,\b \hbox{ and }
\a  \times \b $}
In this section we find precise conditions for the triple 
$\a ,\b \hbox{ and }\a  \times \b $ to be linearly dependant.
Recall that this happens if $\det(\a  \times \b)=0$, ie $\a  \times \b$ is
isotropic. A vector is isotropic if it lies in $X$, 
the right circular cone $x_1^2-x_2^2-x_3^2=0$.

Let 
$$\a\cdot_{\E}\b=a_1b_1+a_2b_2+a_3b_3$$ and 
$$\a\times_{\E}\b=\Biggl|
\matrix{i&j&k\cr a_1&a_2&a_3\cr b_1&b_2&b_3\cr}\Biggr|$$
denote ``euclidean'' scalar and cross product respectively.
This is to distinguish them from the ``hyperbolic'' versions
$$a\cdot_{\H}\b=a_1b_1-a_2b_2-a_3b_3$$ and
$$\a\times_{\H}\b=\Biggl|
\matrix{-i&j&k\cr a_1&a_2&a_3\cr b_1&b_2&b_3\cr}\Biggr|$$
 Let $\rho$ be the involution given by
$\rho(x_1,x_2,x_3)=(-x_1,x_2,x_3)$. Then $X$ is
invariant under $\rho$ and 
$\rho(\a\times_{\H}\b)=\a\times_{\E}\b$.

If $\c$ is isotropic let $\a, \b$ lie in the plane $\c\cdot_{\E}\x=0$.
This plane meets the cone in the generator containing $\rho(\c)$.
So $\a\times_{\E}\b$ is parallel to $\c$ and $\a\times_{\H}\b$ is
isotropic. This means that the triple $\a, \b, \a\times_{\H}\b$ is linearly
dependant. Moreover all examples of such triples are obtained in this way.

From now on $\a\cdot \b$ and $\a\times \b$ will have their original
(hyperbolic) meanings.

\subsection{A Worked Example}

Now we find a generic family of triples containing all the
properties needed. This is summed up by the following theorem
\theorem{Without loss of generality we can assume that $\a$ and $\b$ lie
in the plane $x_1-x_2=0$. The most general examples being
$$\a=\pmatrix{ a_3& 2a_1\cr 0 &-a_3\cr},\quad 
\b=\pmatrix{ b_3& 2b_1\cr 0 &-b_3\cr},\quad 
\a\times\b=\pmatrix{ 0&2(a_3b_1-a_1b_3)\cr 0&0\cr}.$$
$$\hbox{So }b_3\a-a_3\b+\a\times\b=0.$$
}
{\bf Proof}
 Some small lemmas are needed.
\lemma{
  For any traceless matrices $\a, \b, \c$, we have 
$$
    (\a \times\c )\cdot(\b \times\c ) = 
     \det(\c )(\a \cdot\b ) - (\a \cdot\c )(\b \cdot\c )
  $$
}
{\bf Proof} This is just a routine calculation. \qed

We will use conjugation in the group theoretic sense, (ie. $A$ conjugated by
$B$ is $B^{-1}AB$). 
Since the word conjugation is already being used in a rather different sense 
(analogous to complex conjugation) we will 
use the term {\sl group-conjugation}. Note that the set of solutions to the
fundamental equation is invariant under group-conjugation.
 \lemma{ The inner product is invariant under group-conjugation.
This means that for any $2\times2$ matrix $C$ and any tracesless $2\times2$
matrices $\a $ and $\b $ we have 
$$
    C^{-1}\a C\cdot C^{-1}\b C = \a \cdot \b 
  $$
}
{\bf Proof} Tedious but routine calculation using the above lemma. \qed

\lemma{
Any tracesless $2\times2$ matrix is group-conjugate to a matrix of 
the form $a_1\i + a_2\j $. In particular any isotropic tracesless $2\times2$
matrix is group-conjugate to one of the form $x(\i+\j),\ x\in F$}

{\bf Proof} If $a_3 = 0$ then we are already there, so we will assume
otherwise.

Case 1: $a_2^2+a_3^2 \neq 0$.

Consider $C  = (-a_2 + \sqrt{a_2^2+a_3^2})-a_3\i$. Then $C $ is invertible
whenever -
$$
\det(C ) = 2(a_2^2+a_3^2 - a_2\sqrt{a_2^2+a_3^2})=
2(\sqrt{a_2^2+a_3^2})(\sqrt{a_2^2+a_3^2}-a_2) \neq 0 
$$
Hence $C $ is invertible. Moreover
$$
C ^{-1}\a C =  a_1\i + \sqrt{a_2^2+a_3^2}\ \j
$$
so we are done.

Case 2: $a_2^2 + a_3^2 = 0$ 
Again $a_2$ and $a_3$ are not zero and so the underlying field must have a
square root of -1, unique up to multiplication by -1, 
which we will call $I$. In fact $I=\pm a_3/a_2$. Group-conjugating $\a$ by
$1 + I\i$ gives $a_1\i$, as required. \qed

Isotropic vectors are invariant under group conjugation and so by the above
we can
assume an isotropic vector is of the form $x(\i+\j)$. Hence from above $\a$
and $\b$ can be conjugated such that,
$$a_1(\i-\j)  \cdot_{\E} \a = a_1(\i-\j) \cdot_{\E} \b = 0.$$
The most general case then is (up to group-conjugation) 
$\a = a_1\i + a_1\j + a_3\k  \   \b = b_1\i + b_1\j + b_3\k$.
Substituting 
$$\i=\pmatrix{0&1\cr -1&0\cr},\ \j=\pmatrix{0&1\cr 1&0\cr},\ 
\k=\pmatrix{1&0\cr 0&-1\cr}$$
 concludes the proof of Theorem 4.2
\qed

\section{How To Find $A $ and $ B$ given Linearly Dependent $ \a , \b $ and
$\a \times\b $. }
In this section we complete our theoretical solution of $ \a , \b $ to find
$A $ and $ B$. 

Assume that  
$\a ,\b \hbox{ and }\a  \times \b $ are linearly dependent. That is for some
coefficients not all zero there is a linear relationship
$$\lambda_1 \a  + \lambda_2 \b + \lambda_3 \a \times \b=0$$
We are assuming that $A$ and $B$ do not commute, so $\a, \b$ are linearly
independant and hence $\lambda_3 \ne 0$ and we can write
$$\eqalignno{
\lambda_1 \a  + \lambda_2 \b+ \a  \times \b &= 0 &(2)\cr}
$$
where $\lambda_1$ and $\lambda_2 $ are unique. Both  $\lambda_1$ and 
$\lambda_2 $ cannot be zero for then by our earlier discussion $A$ and $B$
would commute.

We can obtain information about $\lambda_1$ and $\lambda_2$ by taking
the cross product of $\a $ with $\a \times \b$
$$\eqalignno{
\a  \times (\a  \times \b ) & = 
         -\lambda_1\a \times\a  - \lambda_2\a \times\b  \cr
          & = \lambda_1 \lambda_2\a  + \lambda_2^2 \b   
\hbox{ on the one hand, and}\cr
          & =   (\b  \cdot \a )\a  - (\a  \cdot \a )\b  
                \hbox{ by standard expansion rules }\cr
} 
$$
Comparing coefficients, we have - 
$$\eqalignno{
 \lambda_1 \lambda_2 &= \a  \cdot \b  ;\cr
 \lambda_2^2 &= -\det(\a ) : \hbox{ and similarly }\cr
 \lambda_1^2 &= -\det(\b )&(3)\cr}
$$

Comparing (2) and (3) with (1) and separating coefficients we get 
the two equations
$$\eqalignno{
\lambda_i[(\det(A)-tr(A))\lambda_1 + b_0\det(A)+2\lambda_1\lambda_2] &= 0
&i=1,\ 2\cr}
$$

Since we cannot have both $\lambda_1$ and $\lambda_2$ = 0, else we get
a commuting solution, we have -
$$\eqalignno{
  b_0 & = {(tr(A)-\det(A)-2\lambda_2)\lambda_1 \over \det(A)} \cr
      & = \lambda_1\left({2 \over a_0+\lambda_2}-1\right)\cr
      & = \sqrt{-\det(\b)}\left({2 \over a_0+\sqrt{-\det(\a)}}-1\right)&(4)\cr
} 
$$

Thus, if $\a $ and $\b $ are such that $\a $, $\b $ and $\a \times\b $
are linearly dependent, then we can pick any $a_0$, not equal to
$-\sqrt{-\det(\a)}$,
then choose $b_0$
according to (4), and we will have $A$ and $B$ that satisfy the
fundamental equation. Note that the roots will have opposite sign.

To be a switch we also require $B$, $A$, $A-I$ and $S$ to be invertible.

Now $A$ is singular if and only if $\det(A) = a_0^2-\lambda_2^2=0$
if and only if $a_0=\pm \lambda_2 $

$(A-1)$ is singular if and only if  $a_0 = 1\pm \lambda_2$

and $B$ is singular  if and only if  
$b_0=\pm\lambda_2$ if and only if  $\lambda_1=0
\hbox{ or } a_0=1-\lambda_2$.

\lemma{If  $A$, $B$ and $(A-1)$ are invertible then so is $S$}

{\bf Proof}
According to \cite{BF}, $S$ is invertible if $\Delta'=C^{-1}D-A^{-1}B$
is invertible. Using $C=A^{-1}B^{-1}A(1-A),\ D=1-A^{-1}B^{-1}AB$
we find that $\Delta'=(1-A)^{-1}A^{-1}B(A-1)$.
\qed

So we require $\lambda_1 \neq 0$ and $a_0 \neq \pm\lambda_2, 1\pm\lambda_2$.

\subsection{The Worked Example(continued)}
 
We have $\lambda_1 = b_3$ and $ \lambda_2 = -a_3$. So
$b_0 = b_3({2 \over a_0-a_3}-1)$
 
hence all mismatching solutions are conjugate to ones of the form -
 $$A = \pmatrix{a_0 + a_3  &  2a_1\cr  0 & a_0 - a_3 \cr}\qquad
 B =  \pmatrix{ {2b_3\over a_0-a_3}   & 2b_1 \cr 0 &  2b_3
\bigl({1\over a_0-a_3}-2\bigr)\cr}
$$ 
The matrices $A$ and $A-I$ will be invertible as long as $a_0 \neq
\pm a_3,\ 2 + a_3\ 1\pm a_3$.

\section{ References}

\refe

\cite{As} Helmer Aslaksen, Quaternionic Determinants, Math. Intel. Vol 18 no.
3 (1996)

\cite{F} Roger Fenn, Quaternion Algebras and Invariants of Virtual Knots
and Links, part I: the Elliptic case, to appear in JKTR

\cite{BF} A. Bartholomew and Roger Fenn. Quaternionic Invariants of Virtual
Knots and Links, to appear in JKTR. Preprint available from\nl
http://www.maths.sussex.ac.uk////Staff/RAF/Maths/Current/Andy/

\cite{BuF} S. Budden and Roger Fenn. The equation $$[b,(a-1)(a,b)]=0$$
and virtual knots and links, Fund Math 184 (2004) pp 19-29.

\cite{C} P. M. Cohn. Algebra vol 3 Wiley 1991.

\cite{FJK} R. Fenn, M. Jordan, L. Kauffman,  Biquandles and 
Virtual Links, Topology and its Applications, 145 (2004) 157-175

\cite{K} L.Kauffman. Virtual Knot Theory, European J. Comb. Vol 20, 
663-690, (1999)

\cite{L} T. Y. Lam. The Algebraic Theory of Quadratic Forms, Benjamin (1973)
\bye

\section{Example}
{$S$: for the Degenerate Case}
\bigskip

$
S=\pmatrix{
a&0&d&0\cr
b&c&e&d(1-c)\cr
(1-a)/d&0&0&0\cr
a(ea-e-bcd)/cd^2(1-c)&1/d&(ea-ec-bcd)/cd(1-c)&0\cr}
$
\bye